\documentclass[12pt]{article}
\usepackage{amssymb,theorem}
mm \textheight=9.0in

\newcounter{gavriil}
\usepackage{amsmath}

\DeclareMathOperator{\diag}{diag}
\DeclareMathOperator{\End}{End}

\DeclareMathOperator{\id}{id}

\DeclareMathOperator{\qtr}{Trace_q}

\DeclareMathOperator{\Span}{Span}

\DeclareMathOperator{\tr}{Trace}

\DeclareMathOperator{\U}{U}

\DeclareMathOperator{\Z}{Z}

\newenvironment{proof}
{\medskip\noindent {\sc Proof:}} {\qed\medskip}

\newcommand{\bm}{\mathbf{m}}

\newcommand{\bn}{\mathbf{n}}

\newcommand{\al}{\alpha}
\newcommand{\be}{\beta}
\newcommand{\De}{\Delta}
\newcommand{\de}{\ensuremath{\delta}}

\newcommand{\Ga}{\Gamma}
\newcommand{\ga}{\gamma}
\newcommand{\ka}{\ensuremath{\varkappa}}

\newcommand{\La}{\ensuremath{\Lambda}}
\newcommand{\la}{\ensuremath{\lambda}}
\newcommand{\Mu}{\mathrm{M}}

\newcommand{\omg}{\omega}

\newcommand{\si}{\sigma}

\newcommand{\ze}{\zeta}

\newcommand{\bla}{\mbox{\boldmath$\lambda$}}
\newcommand{\bmu}{\mbox{\boldmath$\mu$}}
\newcommand{\bnu}{\mbox{\boldmath$\nu$}}

\newcommand{\dl}{{[\![}}
\newcommand{\dr}{{]\!]}}

\newcommand{\dsply}{\displaystyle}

\newcommand{\g}{\mathfrak g}
\newcommand{\gd}{{\mathfrak{g}}^*}

\newcommand{\h}{\mathfrak h}
\newcommand{\hb}{\ensuremath{\hbar}}

\newcommand{\kk}{\mathbb{C}}

\newcommand{\leqs}{\leqslant}

\newcommand{\ot}{\otimes}
\newcommand{\qed}{\ensuremath{\blacksquare}}

\newcommand{\Sg}[1][]{\mathrm{S}^{#1}(\g)}
\newcommand{\tl}{\tilde}
\renewcommand{\to}{\longrightarrow}
\newcommand{\Ug}{{\U}(\g)}
\newcommand{\Ugh}{{\U}_\hb(\g)} 
\newcommand{\Ughr}{{\U}_\hb(\g,r)} 

\newcommand{\w}{\wedge}
\newcommand{\wtl}{\widetilde}

\newcommand{\cA}{\mathcal{A}}
\newcommand{\cB}{\mathcal{B}}

\newcommand{\cF}{\mathcal{F}}

\newcommand{\cH}{\mathcal{H}}

\newcommand{\cJ}{\mathcal{J}}

\newcommand{\cL}{\mathcal{L}}

\newcommand{\cO}{\mathcal{O}}

\newcommand{\cR}{\mathcal{R}}

\newtheorem{thm}{Theorem}
\newtheorem{prop}{Proposition}
\newtheorem{lm}{Lemma}
\newtheorem{cor}{Corollary}

{\theorembodyfont{\rmfamily}\newtheorem{defin}{Definition}}
{\theorembodyfont{\rmfamily}}
{\theorembodyfont{\rmfamily}\newtheorem{rem}{Remark}}

\renewcommand{\ka}{\kappa}

\newcommand{\bmh}{\hat{\bm}}
\newcommand{\bnh}{\hat{\bn}}
\newcommand{\Cart}{\mathfrak{h}} 
\newcommand{\comp}{\circ}   
\newcommand{\dfn}{\em}

\newcommand{\fun}[1][\ola]{\kk[#1]} 
\newcommand{\funh}[1][\ola]{\kk_\hb[#1]} 
\newcommand{\gen}{L}
\newcommand{\gl}{\mathfrak{gl}_n(\kk)}
\newcommand{\gld}{\mathfrak{gl}_n^*(\kk)}
\newcommand{\hcke}{(q-q^{-1})}
\newcommand{\Inv}{Z}
\newcommand{\K}{{\kk[\![\hb]\!]}}

\newcommand{\ola}{O_\la}  
\newcommand{\omu}{O_\mu}  
\newcommand{\oone}{O_1}
\newcommand{\otwo}{O_2}
\newcommand{\pb}{p}  
\newcommand{\Prj}{P}  
\newcommand{\prj}{P}  
\newcommand{\vpi}{\omg}
\newcommand{\vt}{\vartheta}
\newcommand{\vs}{\vspace{.3in}}

\begin{document}

\bibliographystyle{plain}

\title{Quantization of orbit bundles in $\gld$}
\author{A.~Mudrov \\
\small Department of Mathematics, \\
\small University of York, YO10 5DD, UK\\
\small and \\
\small St.-Petersburg Department of Steklov Mathematical Institute,\\
\small Fontanka 27, 191023 St.-Petersburg, Russia
\and V.~Ostapenko\\
\small Department of Mathematics, Bar Ilan University, \\
\small52900 Ramat Gan, Israel }
\date{December 14, 2006}

\maketitle

\begin{abstract}
Let $G$ be the complex general linear group and $\g$ its Lie algebra
equipped with a factorizable Lie bialgebra structure; let $\Ugh$ be
the corresponding quantum group. We construct explicit
$\Ugh$-equivariant quantization of Poisson orbit bundles $\ola
\to \omu$ in $\gd$.
\end{abstract}


e-mail: aim501@york.ac.uk, ostap@math.biu.ac.il

\section*{Introduction}

Let $G$ be a Poisson Lie group and $\g$ the Lie algebra of $G$. Let
$\Ugh$ be the quantization of the universal enveloping algebra $\Ug$
along the corresponding Lie bialgebra structure on $\g$. Consider
two Poisson homogeneous $G$-manifolds, $M:=G/G_M$ and $N:=G/G_N$
such that the stabilizer $G_M$ is a subgroup in $G_N$. Then there
exists a natural projection $M \to N$ of $G$-spaces, and it makes
$M$ into a $G$-bundle over $N$ with the fiber $G_N/G_M$. Suppose
this projection is a Poisson map. It is natural to consider the
problem of equivariant quantization of such a bundle. By this we
understand a $\Ugh$-equivariant quantization of the function algebras
together with the co-projection $F(M)\hookleftarrow F(N)$, to a
morphism $F_\hb (M)\hookleftarrow F_\hb (N)$ of $\Ugh$-algebras.

In the present paper, we quantize orbit bundles for the case when
$G=GL_n(\kk)$, and $\g$ is equipped with a factorizable Lie
bialgebra structure. We assume the stabilizers $G_M$
and $G_N$ to be Levi subgroups of $G$.
Specifically for the $GL_n(\kk)$-case, those are precisely
reductive subgroups of maximal rank.
Then the $G$-varieties $M$ and $N$ can be
realized as semisimple coadjoint orbits  $\oone,\otwo\subset\gd$.

The Poisson structure on $\cO_i$ is obtained by restriction
from a Poisson structure on $\gd$. The latter is a linear combination of
the $G$-invariant
Kostant-Kirillov-Lie-Souriau bracket (KKLS) and the
Semenov-Tian-Shansky bracket (STS).
In our case, $G=GL_n(\kk)$, they
 are compatible, i.e. the Schouten bracket
between the two is equal to zero.
The Poisson bracket on  $\cO_i$ is not $G$-invariant, however, it makes
$\cO_i$ Poisson-Lie manifolds over the Poisson-Lie group $G$.
Moreover, it is the only such bracket on  $\cO_i$ obtained by restriction from $\gd$.

Explicit quantization of semisimple orbits has been constructed in
\cite{DoninJMudrovAI:ExplicitQuantiz} for the special case of the
standard, or Drinfeld-Jimbo, quantum group $\Ugh$.
In the present paper, we extend that quantization  to the case of any
factorizable Lie bialgebra structure on $\g$ and the
corresponding quantum group.
We also describe all semisimple Poisson-Lie orbit bundles
$\oone\to\otwo$ in $\gd$.
In particualar, we show that $\otwo$ is necessarily
symmetric.
We explicitly construct a $\Ugh$-equivariant quantization
of the projection map $\prj$ for all orbit bundles.

The paper is organized as follows.

\noindent
In Section~1 we recall some
definitions concerning equivariant quantization.

\noindent
In Section~2 we establish necessary and sufficient conditions for
an orbit bundle to be Poisson.

\noindent
In Section~\ref{sect:REAlgebras} we study the behavior of algebras
defined by a modified (quadratic-linear) Reflection Equation
under twist of quantum groups.

\noindent
In Section~\ref{sect:DMQuantization} we use results of the previous
section to extend the double quantization of orbits
\cite{DoninJMudrovAI:ExplicitQuantiz} to the case of the quantum
group defined by an arbitrary factorizable classical r-matrix.

\noindent
In Section~\ref{sect:QnOfBundles} we prove that any Poisson orbit
bundle admits a $\Ugh$-equivariant quantization, so the conditions
of Section \ref{sect:Poisson-LieOrB} are also sufficient.
We give an explicit formula for the quantized bundle map.

\noindent
There is an Appendix at
the end of the paper where we study certain properties of the q-trace functions.

\subsubsection*{Acknowledgements}
We thank Steven Shnider whose valuable remarks led to a
substantial improvement of the text.
This work was supported by the EPSRC grant C511166
 and by the RFBR grant 06-01-00451.

\section{Generalities on equivariant quantization}

\subsection{Deformation quantization of Poisson varieties}

Let $M$ be a variety with a Poisson bracket $\pb$ and $A=\fun[M]$ be the
algebra of polynomial functions $M\to\kk$.
Recall the following definition (see e.g. \cite{DrinfeldW:Doklad}):
\begin{defin}
\label{def:QnOfPoissManifs} An algebra $(A_\hb,\star)$ over the ring
$\K$ of formal power series is called {\dfn quantization} of
$(A,\pb)$ if:
\begin{enumerate}
\item[ (i)]
$A_\hb$ is a free $\K$-module;
\item[ (ii)]
As a $\kk$-algebra, the quotient $A_\hb/\hb A_\hb$ is isomorphic to
$A$;
\item[ (iii)]
If $a,b\in A$ then $\dsply\frac{a\star b-b\star
a}{\hb}\equiv\pb(a,b)$ modulo $\hb$.
\end{enumerate}
The Poisson bracket $\pb$ is called {\dfn the infinitesimal} of
$(A_\hb,\star)$.
\end{defin}
\begin{rem}\label{nb:Pb}
The deformed multiplication is expanded as an $\hbar$-series:
$\dsply a\star b=\sum_{k=0}^\infty m_k(a\ot b)\hb^k$
for $a,b\in A\subset A_\hbar$.
Therefore one has  $\pb(a,b)=m_1(a\ot b) - m_1(b\ot a)$.
\end{rem}

Let $(M,\pb_M)$ and $(N,\pb_N)$ be two Poisson varieties, $A_\hb$
and $B_\hb$ some quantizations of the function algebras $A=\fun[M]$
and $B=\fun[N]$ respectively, $f\colon B\to A$ a morphism of Poisson
algebras.
\begin{defin}
\label{def:QnOfPoissMorphs}
A homomorphism $f_\hb\colon B_\hb\to A_\hb$ of
$\K$-algebras is called a {\dfn quantization of the map} $f$ if the
induced  morphism $f_0\colon B_\hb/\hb B_\hb\to A_\hb/\hb
A_\hb$ of $\kk$-algebras coincides with $f$.
\end{defin}
\begin{prop}
Suppose there exists a quantization of $f\colon B\to A$.
Then $f$ is a Poisson map.
\end{prop}
\begin{proof}
Denote by $m_{A_\hb}$ the multiplication in $A_\hb$ and by $m_A$ the
undeformed multiplication in $A$. The map $f_\hb$ is supposed to
be an algebra homomorphism, i.e. $f_\hb\comp
m_{B_\hb}=m_{A_\hb}\comp\left(f_\hb\ot f_\hb\right)$. Consider the
infinitesimal part of this equality:
\begin{equation*}
f\comp m_{B_\hb,1}+f_1\comp m_B= m_A\comp\left(f_1\ot f\right)+
m_A\comp\left(f\ot f_1\right)+ m_{A_\hb,1}\comp\left(f\ot f\right).
\end{equation*}
Applying this equality first to $a\ot b\in A\ot A$, then to $b\ot a$
and taking the difference, one obtains $f\bigl(\pb_M(a,b) \bigr) =
\pb_N\bigl(f(a),f(b)\bigr)$ where $\pb_M(a,b)=m_{A_\hb,1}(a,b)-m_{A_\hb,1}(b,a)$ and
$\pb_N(a,b)=m_{B_\hb,1}(a,b)-m_{B_\hb,1}(b,a)$.
\end{proof}

\subsection{Quantization of $G$-varieties}
Consider  a simple complex algebraic group  $G$ and its Lie algebra  $\g$.
Suppose $G$ is a Poisson group; then $\g$  is equipped with a quasitriangular
Lie bialgebra structure. Denote by $r\in \wedge^2\g$
the corresponding classical r-matrix.
We consider only the factorizable case, when $r$ satisfies the {\em modified} classical Yang-Baxter equation.
By $\Ughr$ we denote the quantization of the universal enveloping algebra $\Ug$ along
$r$.

Consider the variety $M=K\backslash G$ where $K\subset G$ is a reductive subgroup in $G$ of maximal rank;
let $\mathfrak{k}$ denote its Lie algebra (we prefer to work with right coset spaces, so that the right $G$-action induces
a left action on functions).
Suppose $M$ is a Poisson $G$-variety, that is to say,
the action $G\times M\to M$ is Poisson. Set $A=\kk[M]$  and let $A_\hb$ be its quantization.
We expect the deformed
multiplication in $A_\hb$ to be equivariant with respect to an
action of $\Ughr$. In other words, this multiplication should obey the
"Leibniz rule"
$$
x.(a\star b)=(x^{(1)}.a)\star (x^{(2)}.b)
$$
for all $a,b\in A$ and $x\in\Ughr$.
We use the  standard Sweedler notation $x^{(1)}\ot x^{(2)}$
for the coproduct $\De(x)$.

The infinitesimal of a
$\Ughr$-equivariant quantization $A_\hb$
is always of the form (see \cite{DoninGurevicShnider:DoubleQuantiz})
\begin{equation}\label{eq:UhEquivariantPb}
\pb=\overleftarrow{r}+\overrightarrow{s}.
\end{equation}
Here $\overleftarrow{r}$ denotes the bivector field on
$M$ generated by $r\in \wedge^2\g$ via the action of $G$,
and $\overrightarrow{s}$ denotes the  invariant bivector field on
$M$ generated by an element $s\in\left(\wedge^2\g/\mathfrak{k}\right)^\mathfrak{k}$.
The latter
should satisfy certain conditions, and the bivector field
$\overrightarrow{s}$  is called {\dfn quasi-Poisson structure} on $M$.
For simplicity, we call the generator $s$  a quasi-Poisson structure as well.

Recall the constructions of the bivector fields $\overleftarrow{r}$
and $\overrightarrow{s}$.
For $r\in\wedge^2\g$ and $g\in G$, consider the bivector field $(R_g)_*r$
on $G$,
where $R_g:G\to G$ is the right translation $x\mapsto xg$.
This bivector field is left $G$-invariant, so it is left $K$-invariant.
Hence it is projectable to $M=K\backslash G$, and we denote the projection by
$\overleftarrow{r}$.
Note that $\overleftarrow{r}$ is not $G$-invariant.

To describe the bivector field $\overrightarrow{s}$,
lift $s$ from $\wedge^2\left(\g/\mathfrak{k}\right)^\mathfrak{k}$
to  $\left(\wedge^2\g\right)^\mathfrak{k}$ and consider
a left $G$ invariant bivector field $(L_g)_*s$,
where $L_g:G\to G$ is the left translation $x\mapsto gx$.
It is also left $K$-invariant, hence it is projectable to $M=K\backslash G$.
We denote the projection by $\overrightarrow{s}$.
Any $G$ invariant bivector field on $M$ is obtained in this way.

Left and right invariant vector
fields on $G$  commute with each other, hence the Schouten bracket
$[\overleftarrow{r},\overrightarrow{s}]$ vanishes.
This implies that  $p=\overleftarrow{r}+\overrightarrow{s}$ is
a Poisson bracket if and only if  $[\overleftarrow{r},\overleftarrow{r}]=-[\overrightarrow{s},\overrightarrow{s}]$.

Recall that two Poisson brackets are called compatible if their any linear combination
is again a Poisson bracket. Suppose that the new bracket makes the variety $M$ Poisson over the Poisson
group $G$. Then the formula (\ref{eq:UhEquivariantPb}) suggests that
a Poisson bracket $\ka$ on $M$ is compatible with $p$ if and only if it is
$G$-invariant and $[\overrightarrow{s},\ka]=0$.
Next we recall the notion of 2-parameter, or double quantization (see,
for example, \cite{DoninJ:DoubleQuantizationOfLieGroups}).

\begin{defin}
\label{def:DoubleQnOfPoissManifs} Suppose that the commutative
algebra $A$ is endowed with two compatible Poisson brackets, $\pb$
and $\ka$, such that $\ka$ is $G$-invariant. An algebra
$(A_{\hb,t},\star)$ over the ring $\kk\dl\hb,t\dr$ of formal power
series in two variables is called {\dfn equivariant quantization} of
$(A,\pb,\ka)$ if
\begin{enumerate}
\item[ (i)]
$A_{\hb,t}$ is a free module over $\kk\dl\hb,t\dr$;
\item[ (ii)]
The first order term of the deformed multiplication $\star$ is
$\hb\pb+t\ka$;
\item[ (iii)]
The algebra $A_{\hb,t}$ is $\Ughr$-equivariant;
\item[ (iv)]
The quotient $A_{\hb,t}/\hb A_{\hb,t}$ is a $G$-equivariant
(one-parameter) quantization of $(A,\ka)$.
\end{enumerate}
\end{defin}

\section{Poisson-Lie orbit bundles}
\label{sect:Poisson-LieOrB}
\subsection{General remarks on coadjoint orbits in $\gld$}
\label{sect:GenRemarks}

Fix $G=GL_n(\kk)$ and put $\g=\gl$. Choose the
algebra of diagonal matrices as a Cartan subalgebra
$\Cart\subset \g$.
The non-degenerate $G$-invariant bilinear
form $\tr(XY)$ on $\g$ allows us to think of $\Cart^*$
as a subset in $\g^*$, the dual vector space for $\g$.
Denote by $\ola$ the coadjoint orbit of a
semisimple element $\la\in\Cart^*\subset\gd$.
As a $G$-variety,
$\ola$ is isomorphic to $G^\la\backslash G$ where $G^\la$ is a Levi
subgroup of $G$. That is, the Lie algebra $\g^\la$ of $G^\la$
is a Levi subalgebra of $\g$ containing
$\Cart$.

The trace bilinear form on $\g$ establishes a canonical
isomorphism between coadjoint orbits in $\gd$ and adjoint orbits
in $\g$.
We will use this isomorphism without further noticing.
In the same way, we identify an element $\la\in\Cart^*$
with the corresponding diagonal matrix in $\Cart$.

Two diagonal matrices with different order of the entries belong
to the same $G$-orbit.
Hence we can
choose a representative of the orbit in which all equal entries are
grouped up together. In other words, we can think that
$\la=\diag(\La_1,\ldots,\La_l)$ where $\La_i$ is the scalar
$n_i\times n_i$ matrix with $\la_i$ on its diagonal.
In particular, the orbit $\ola$ is determined by a pair
$(\bla,\bn)$ where $\bla=(\la_1,\ldots,\la_l)$ is the row of pairwise
distinct eigenvalues of $\la$, and
$\bn=(n_1,\ldots,n_l)$ and $\la_1,\dots,\la_l$ is the row of their
multiplicities.
Note that the orbit $\ola$ does not depend on
simultaneous permutations of the entries of $\bla$ and $\bn$.

\subsection{The related Poisson structures on $\ola$}

The admissible quasi-Poisson bracket on a semisimple orbit has the form
(see  \cite{DoninJ:QuanGMans},
\cite{DoninJOstapenkoV:EquivQn},
\cite{DoninGurevicShnider:DoubleQuantiz} and
\cite{KarolinskijE:ClassPoissHomSpace} for details):
$$
\sum_{1\leqs i<j\leqs l}c_{ij}\xi_{ij},
$$
where $c_{ij}$ are some coefficients depending on the eigenvalues of
$\la$ (see below for explicit formulas),
and $\xi_{ij}$ are defined as follows:
\begin{equation}\label{eq:BasisBivector}
\xi_{ij}=\sum_{s,t} E_{st}\w E_{ts} \mod \g\wedge\g^\la.
\end{equation}
Here $E_{st}$ is the $(s,t)$-th matrix unit, $l$ is the number of
different eigenvalues of $\la$, and the sum is taken over
$n_1+n_2+\dots+n_{i-1}< s \leqs n_1+n_2+\dots+n_i$,
$n_1+n_2+\dots+n_{j-1}< t \leqs n_1+n_2+\dots+n_j$, where $n_i$
denotes the multiplicity of the eigenvalue $\la_i$.
\subsubsection*{The KKSL Poisson bracket $\ka_\la$}

The Kirillov-Kostant-Lie-Souriau bra\-cket $\ka_\la$ is induced on
$\ola$ from the Lie structure on $\g$.
For a semisimple (co)adjoint $GL_n(\kk)$-orbit $\ola$ it is given
by the following expression:
\begin{equation*}
\ka_\la=\sum_{1\leqs i<j\leqs l}\frac{1}{\la_i-\la_j}\ \xi_{ij}.
\end{equation*}
This is a $G$-invariant non-degenerate Poisson bracket.
\subsubsection*{The quasi-Poisson bracket $s^0_\la$.}
Consider the bivector field on $\ola$ restricted
from the Semenov-Tian-Shansky (STS)
bracket on $\End(\kk^n)$, see \cite{STS:QDualityQDouble}.
Its quasi-Poisson part is generated by an element $s_\la^0\in\wedge^2(\g/\g^\la)^{\g^\la}$.
Specifically for the
case $G=GL_n(\kk)$, it takes the form:
$$
s_\la^0=\sum_{1\leqs i<j\leqs l} \frac{\la_i+\la_j}{\la_i-\la_j}\
\xi_{ij},
$$
see  \cite{DoninJMudrovAI:MethodQuantChar}.

\subsubsection*{The brackets admitting $\Ughr$-equivariant
quantization}

In our case $G=GL_n(\kk)$, the Schouten bracket between $s_\la^0$
and $\ka_\la$ vanishes.
Thus, by adding a multiple of the KKLS
bracket to that generated by $s_\la^0$, one obtains the general form
for a quasi-Poisson bracket on $\ola$ admitting a $\Ughr$-equivariant
quantization:
\begin{equation}\label{eq:sla}
s_\la=s_\la^0+a \ka_\la= \sum_{1\leqs i<j\leqs
l}\frac{\la_i+\la_j+a}{\la_i-\la_j}\ \xi_{ij}, \mbox{ where }
a\in\kk.
\end{equation}
Recall once again that we consider only those Poisson structures
that are restricted from $\gd$.
\subsection{The structure of Poisson orbit bundles}

Here  we give the necessary and sufficient conditions for
an orbit map $\prj\colon \left(\ola,s_\la\right)\to\left(\omu,s_\mu\right)$
to be Poisson.
Recall that we do not distinguish between an element $\la\in\gd$ and
the corresponding diagonal matrix $\diag(\La_1,\ldots,\La_l)$.
\begin{prop}\label{prop:P}
Any semisimple orbit bundle $\ola\to\omu$ is of the form
$X\mapsto\Prj(X)$ where
$\Prj$ is a polynomial in one variable with complex coefficients.
\end{prop}
\begin{proof}
If an equivariant map $\left(\ola,s_\la\right)\to\left(\omu,s_\mu\right)$
brings $\la$ to $\mu$, then the isotropy group $G^\la$ is
a subgroup in $G^\mu$; this gives an inclusion $\g^\la\subset \g^\mu$ of their Lie algebras.

Consider $\g=\gl$ as an associative algebra, $\g\cong\End(\kk^n)$,
and denote by $\Z(\g^\la)\subset\g$ the centralizer of
$\g^\la$.
Since $\g^\la$ is a Levi subalgebra,
$\Z(\g^\la)$ is a semisimple  commutative associative
algebra generated by $\la$ and by the unit matrix.
For any polynomial
$P$ in one variable the mapping $X\mapsto P(X)$, $X\in \g$, is
$G$-equivariant.
Hence it suffices to check that $\mu = \Prj(\la)$
for some polynomial $\Prj$.
The inclusion $\g^\la\subset \g^\mu$ implies the inclusion $Z(\g^\mu)\subset Z(\g^\la)$.
Therefore the matrix $\mu$ is a polynomial in
$\la$.
\end{proof}
\begin{rem}
We will use the same symbol $\Prj$ for both the orbit map
$\ola\to\omu$ and for the corresponding polynomial.
\end{rem}
\subsubsection*{The graph of an orbit bundle}

It is convenient to use
the following graphical presentation of an orbit bundle $\Prj\colon\omu\to\ola$.
The orbit $\ola$ has a representative in  the form of diagonal matrix
$\la=\diag(\La_1,\ldots,\La_l)$.
Set $\mu=\Prj(\la)$, then $\mu$ also has the form
$\mu=\diag(\Mu_1,\ldots,\Mu_m)$, where $\Mu_j$ denotes the scalar
block corresponding to the eigenvalue $\mu_j$ of $\mu$.

Using this, denote by $\Ga_\Prj$ the bipart type
graph whose upper nodes have labels $1,\dots,l$ corresponding to the
blocks $\La_1,\dots,\La_l$, and the lower nodes are labeled by
$1,\dots,m$ corresponding to the blocks $\Mu_1,\dots,\Mu_m$.
The $i$-th node of the upper part of $\Ga_\Prj$ is connected to the
$\al$-th node of the lower part if and only if
$\Prj(\la_i)=\mu_\al$.
Note that each upper node has exactly one edge.
Since the map $\Prj$ is surjective, each
lower node is connected to some upper node.
Therefore the graph $\Ga_\Prj$ is a disjoint union of trees of
the form:
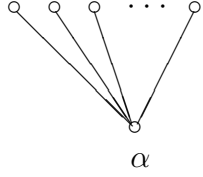
\begin{center}\setlength{\unitlength}{.7in}
\begin{picture}(0,-1)
  \multiput(-.9,0)(.3,0){3}{\circle{.08}}
  \put(.45,0){\circle{.08}}
  \put(0,-.9){\circle{.08}}
  \put(-.88,-.03){\line(1,-1){.863}}
  \put(-.59,-.03){\line(2,-3){.565}}
  \put(-.3,-.03){\line(1,-3){.279}}
  \put(.445,-.03){\line(-1,-2){.423}}
  \put(-.03,-1.2){$\al$}

  \multiput(-.06,-.05)(.12,0){3}{$\cdot$}
\end{picture}
\end{center}\vspace{1.0in}
Each tree is a connected component of $\Ga_\Prj$
and it is labeled by the blocks of $\mu$.
The graph  $\Ga_\Prj$ gives a complete description of
the map $\Prj$.
\subsubsection*{The graph $\Ga_\Prj$ for $\Prj$ Poisson}

Both $\ola$ and $\omu$ are endowed with Poisson
structures (\ref{eq:sla}).
The tangent space of $\ola$ at the point $\la$ is isomorphic as
a vector space to the quotient $\g/\g^\la$.
Recall that $\g^\la\subset \g^\mu$.
The tangent map $\prj_*:\g/\g^\la\to\g/\g^\mu$ of $\prj$ is given by
the formula:
$$
\prj_*(X)=\left\{
\begin{array}{l}
0, \mathrm{\ if\ } X\in\g^\mu/\g^\la,\\
\ \\
X \mathrm{\ otherwise.}
\end{array}
\right.
$$
An element $X\in\g\setminus\Cart$ cannot generate a vector field on
$\omu$ since it is not even $\Cart$-invariant.
However, the element like $\xi_{ij}$ (see formula
(\ref{eq:BasisBivector})) does generate a bivector field on $\omu$.
The map $\prj$ is  Poisson if and only if
$\prj_*(s_\la)=s_\mu$.
The tangent map $\prj_*$ is determined by
its values $\prj_*(\xi_{ij})=\xi_{\al\be}$ where $i,j$ run over
the upper nodes of the graph $\Ga_\Prj$ while $\al,\be$ run
over the lower nodes of $\Ga_\Prj$.
\begin{lm}\label{lm:GraphOfPoissBundle}
Let $\Prj\colon \ola\to\omu$ be a Poisson  orbit bundle w.r.t.
the Poisson structure on $\ola$ determined by $s_\la=\sum
c_{ij}\xi_{ij}$ and some Poisson structure on $\omu$.
Suppose that $i\leqs j<p\leqs s$ are such {\em upper} nodes that $i$
and $j$ belong to the same connected component $\al$, and $p$, $s$ also
belong to the same connected component $\be$ of $\Ga_\Prj$.
Then $c_{ip}=c_{js}$.
\end{lm}
\begin{proof}
If $\Prj$ is Poisson then $s_\la$ is projectable under $\prj_*$.
Thus if $\prj_*(\xi_{ip})$ and $\prj_*(\xi_{js})$
enter the same basis bivector $\xi_{\al\be}$, then $c_{ip}=c_{js}$.
\end{proof}
\subsubsection*{Classification of Poisson bundles}

We call a coadjoint orbit $\omu$ {\em symmetric} if the
corresponding matrix $\Mu$ has exactly two different eigenvalues.
\begin{thm}\label{th:MainPoiss}
A $GL_n(\kk)$-equivariant map $\prj\colon \oone\to\otwo$ is Poisson if and
only if the following three conditions are satisfied:
\begin{enumerate}
\item[{\em (a)}]
The orbit $\otwo$ is symmetric;
\item[{\em (b)}]
There exist $\la\in\oone$ and $\mu\in\otwo$ such that
$\prj(\la)=\mu$ and the multiplicity $n_1$ of the eigenvalue $\la_1$
is equal to the multiplicity $m_1$ of the eigenvalue $\mu_1$;
\item[{\em (c)}]
The Poisson structures on $\oone$ and $\otwo$ are defined by $s_\la
= s_\la^0 - 2\la_1 \ka_\la$ and $s_\mu = s_\mu^0 - 2\mu_1 \ka_\mu$
respectively.
\end{enumerate}
\end{thm}
\begin{proof}
Show first that if $\Prj\colon \oone\to\otwo$ is a Poisson map,
then the orbit $\otwo$ is symmetric, i.e. the graph $\Ga_\Prj$
consists of exactly two connected components.
Indeed, suppose that $\Ga_\Prj$ has the
form \vs
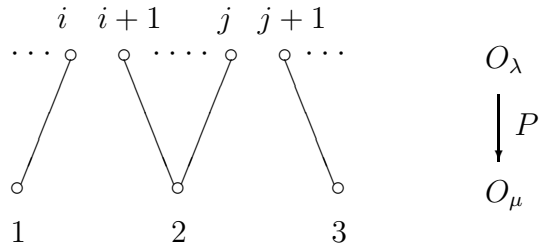
\begin{center}\setlength{\unitlength}{.7in}
\begin{picture}(0,-2)
  \put(-1.6,0){\circle{.08}}
  \put(-1.2,0){\circle{.08}}
  \put(-.4,0){\circle{.08}}
  \put(0,0){\circle{.08}}

  \put(-2.0,-1){\circle{.08}}
  \put(-.8,-1){\circle{.08}}
  \put(.4,-1){\circle{.08}}

  \put(-1.6,-.04){\line(-2,-5){.37}}
  \put(-1.2,-.04){\line(2,-5){.37}}
  \put(-.4,-.04){\line(-2,-5){.37}}
  \put(0,-.04){\line(2,-5){.37}}

  \put(-1.7,.2){$i$}
  \put(-1.4,.2){$i+1$}
  \put(-.5,.2){$j$}
  \put(-.2,.2){$j+1$}

  \put(-2.05,-1.4){$1$}
  \put(-0.85,-1.4){$2$}
  \put(.35,-1.4){$3$}

  \multiput(-2.05,-.05)(.12,0){3}{$\cdot$}
  \multiput(-.99,-.05)(.12,0){4}{$\cdot$}
  \multiput(.15,-.05)(.12,0){3}{$\cdot$}

  \put(1.5,-.1){$\ola$}
  \put(1.5,-1.1){$\omu$}
  \put(1.6,-.3){\vector(0,-1){.5}}
  \put(1.72,-.6){$\Prj$}
\end{picture}
\end{center}\vs\vs\vs
Since $\Prj$ is Poisson,  Lemma~\ref{th:MainPoiss} leads to the
following system of conditions:
$$
\left\{
\begin{array}{l}
c_{i,i+1}=c_{i,j}\\
c_{i+1,j+1}=c_{j,j+1},
\end{array}
\right.
$$
where $\dsply c_{ij}=\frac{\la_i+\la_j+a}{\la_i-\la_j}$, see formula
(\ref{eq:sla}). Thus we should solve the following overdefined
system of linear equations in $a$:
$$
\left\{
\begin{array}{l}
\dsply \frac{\la_i+\la_{i+1}+a}{\la_i-\la_{i+1}}=
\frac{\la_{i}+\la_{j}+a}{\la_{i}-\la_{j}} \\
\   \\
\dsply \frac{\la_{i+1}+\la_{j+1}+a}{\la_{i+1}-\la_{j+1}}=
\frac{\la_{j}+\la_{j+1}+a}{\la_{j}-\la_{j+1}}
\end{array}
\right.
$$
This system is inconsistent, thus our hypothesis is wrong, and
$\Ga_\Prj$ consists of exactly two connected trees.
This means that
$\mu$ has exactly two different eigenvalues:
$\mu=\diag(\Mu_1,\Mu_2)$, i.e. it is symmetric.

We now show that if $\Prj$ is Poisson, then either
$\mathrm{mult}\>\la_1=\mathrm{mult}\>\mu_1$ or $\mathrm{mult}\>\la_l=\mathrm{mult}\>\mu_2$,
where $\mathrm{mult}$ means the multiplicity of the corresponding eigenvalue.
In other words, we need to prove that
if $\Prj$ is Poisson, then one of the two connected components of
the graph $\Ga_\Prj$ contains precisely one edge.
Suppose, to the contrary, that $\Ga_\Prj$ is of the form:
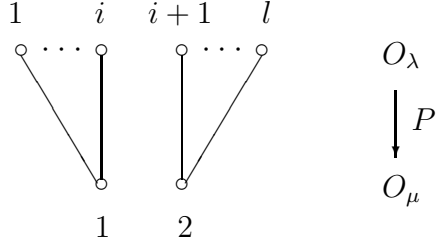
\vs
\begin{center}\setlength{\unitlength}{.7in}
\begin{picture}(0,-1)
  \put(-1.2,0){\circle{.08}}
  \put(-.6,0){\circle{.08}}
  \put(0,0){\circle{.08}}
  \put(.6,0){\circle{.08}}

  \put(-.6,-1){\circle{.08}}
  \put(0,-1){\circle{.08}}

  \put(-1.2,-.04){\line(3,-5){.57}}
  \put(-.6,-.04){\line(0,-1){.93}}
  \put(0,-.04){\line(0,-1){.93}}
  \put(.6,-.04){\line(-3,-5){.57}}

  \put(-1.3,.2){$1$}
  \put(-.65,.2){$i$}
  \put(-.25,.2){$i+1$}
  \put(.6,.2){$l$}

  \put(-.65,-1.4){$1$}
  \put(-.03,-1.4){$2$}

  \multiput(-1.05,-.05)(.12,0){3}{$\cdot$}
  \multiput(.15,-.05)(.12,0){3}{$\cdot$}

  \put(1.5,-.1){$\ola$}
  \put(1.5,-1.1){$\omu$}
  \put(1.6,-.3){\vector(0,-1){.5}}
  \put(1.72,-.6){$\Prj$}
\end{picture}
\end{center}\vs\vs\vs
Again, by Lemma~\ref{th:MainPoiss} the following system of
conditions should be satisfied: $c_{1,i+1}=c_{i,i+1}=c_{i,l}$.
This is a system of equations in $a$:
$$
\frac{\la_{1}+\la_{i+1}+a}{\la_{1}-\la_{i+1}}=
\frac{\la_i+\la_{i+1}+a}{\la_i-\la_{i+1}}=
\frac{\la_{i}+\la_{l}+a}{\la_{i}-\la_{l}},
$$
which obviously has no solution.
Thus one of the two connected
components of $\Ga_\Prj$ contains one edge.
By appropriate choice of
the representatives $\la\in\ola$ and $\mu\in\omu$, one can assume
that to be the left component:
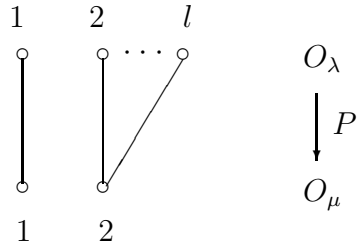
\vs
\begin{center}\setlength{\unitlength}{.7in}
\begin{picture}(0,-1)
  \put(-.6,0){\circle{.08}}
  \put(0,0){\circle{.08}}
  \put(.6,0){\circle{.08}}

  \put(-.6,-1){\circle{.08}}
  \put(0,-1){\circle{.08}}

  \put(-.6,-.04){\line(0,-1){.93}}
  \put(0,-.04){\line(0,-1){.93}}
  \put(.6,-.04){\line(-3,-5){.57}}

  \put(-.7,.2){$1$}
  \put(-.1,.2){$2$}
  \put(.6,.2){$l$}

  \put(-.65,-1.4){$1$}
  \put(-.03,-1.4){$2$}

  \multiput(.15,-.05)(.12,0){3}{$\cdot$}

  \put(1.5,-.1){$\ola$}
  \put(1.5,-1.1){$\omu$}
  \put(1.6,-.3){\vector(0,-1){.5}}
  \put(1.72,-.6){$\Prj$}
\end{picture}
\end{center}\vs\vs\vs
i.e. that the multiplicities of $\la_1$ and $\mu_1$
coincide.

Lemma \ref{lm:GraphOfPoissBundle} applied to the above graph gives rise
to the overdefined system of $l-1$ linear equations
$c_{12}=c_{13}=\ldots=c_{1l}$. This system is
consistent, and the solution is $a=-2\la_1$. Hence the
quasi-Poisson brackets on $O_i$ are as in (c).

Now we prove the sufficiency of the conditions (a)--(c).
Indeed, consider an invariant mapping  $\Prj\colon \oone\to\otwo$
corresponding to the above graph.
Suppose that $\Prj(\la)=\mu$. Then for $s_\la =
s_\la^0 - 2\la_1 \ka_\la$ and $s_\mu = s_\mu^0 - 2\mu_1 \ka_\mu$ one
has $\Prj_*(s_\la)=s_\mu$, i.e. $\Prj$ is a Poisson map with respect
to those brackets.
\end{proof}

\subsection{Explicit formula for the map $\prj$}

Denote by $\la$ and $\mu$ elements of
$\Cart$ satisfying the conditions of Theorem~\ref{th:MainPoiss} and by
$\ola$ and $\omu$ their adjoint orbits.
For the purpose of quantization, we need an explicit expression for
$\prj$.
Recall that $\la_1,\la_2,\dots,\la_l$ denote all the distinct
eigenvalues of $\la\in\g$, so we need to find a polynomial
$\Prj(x)$ of degree $l-1$ in one variable with complex coefficients
such that
$\Prj(\la_1)=\mu_1$ and $\Prj(\la_2)=\ldots=\Prj(\la_l)=\mu_2$.
These $l$ equations determine $\Prj$ uniquely:
\begin{equation}\label{eq:DefOfP}
\Prj(x)=(\mu_1-\mu_2)\prod_{i=2}^l\frac{x-\la_i}{\la_1-\la_i} +
\mu_2.
\end{equation}

\section{Reflection equation algebras}
\label{sect:REAlgebras}

An equivariant two-parameter quantization of coadjoint orbits is
constructed in \cite{DoninJMudrovAI:ExplicitQuantiz} for the
special case of the standard, or Drinfeld-Jimbo, quantum group
$\Ugh$.
In this section we extend the quantization of \cite{DoninJMudrovAI:ExplicitQuantiz}  for an arbitrary
quantum group $\Ughr$, not necessarily the standard.
Note that
possible factorizable Lie bialgebra structures on $\g$ are
parameterized by  Belavin-Drinfeld triples and a subspace in $\h$,
 \cite{BelawinDrinfeld:CYBE}.
Quantization of the universal enveloping algebra along any Lie
bialgebra structure has been constructed in
\cite{EtinhoffKazdan:QuantBialg}.

To describe quantized coadjoint orbits explicitly, we need some
facts about the so called {\em modified reflection equation} (mRE)
algebra.
It is a two parameter quantization of the polynomial ring
on the vector space of $n\times n$-matrices.
The quantized orbits will be
presented as quotients of the mRE algebra by certain ideals which
are deformations of the classical ideals of the orbits.
In the present section we study how mRE algebras transform under
twist of quantum groups.
\subsection{The definition and basic facts}
\label{subsection:mREA}

Let $\Ughr$ be a quantization of the universal enveloping algebra of
$\g=\gl$ along a classical r-matrix $r$.
Let $R\in \End(\kk^n\ot
\kk^n)\dl\hb\dr$ be the image of its universal R-matrix in the basic
representation in $\kk^n\dl\hb\dr$.
Denote by $S:=\si R\in
\End(\kk^n\ot\kk^n)\dl\hbar\dr$ the quantum permutation, where $\si$
designates the usual flip $\kk^n\ot\kk^n\to \kk^n\ot\kk^n$, $u\ot
v\mapsto v\ot u$.
\begin{defin}
The mRE algebra $\cL$ is an associative unital
algebra over the ring $\kk\dl\hbar\dr[t]$ generated by the entries of
the matrix $\gen=\left(L_{ij}\right)_{i,j=1}^n$
modulo the relations
\begin{equation}\label{eq:mRE}
[S\gen_2S,\gen_2]=-qt[S,\gen_2] \in \End(\kk^n\ot
\kk^n)\ot\cL,
\end{equation}
where $\gen_2:=1\ot \gen$ and $q:=e^\hbar$.
\end{defin}
The action of $\Ughr$ on the algebra $\cL$ is given by the formula
\begin{equation}\label{eq:action}
x\triangleright L= \rho\bigl(\ga(x^{(1)})\bigr)L\rho(x^{(2)}),
\end{equation}
where $\ga$ is the antipode of $\Ughr$ and $\rho$ is the representation
$\Ughr\to\End(\kk^n)\dl\hb\dr$.
The algebra $\cL$ is a $\Ughr$-equivariant
quantization of the polynomial ring $\kk[\End(\kk^n)]$
with $\hb s^0_\la-t\kappa_\la$
being the linear term of the deformed multiplication,
\cite{DoninJMudrovAI:ExplicitQuantiz}.

\begin{rem}
The relations (\ref{eq:mRE}) are called the
{\em modified reflection equation} (mRE).
These relations become quadratic when $t=0$.
The corresponding quotient $\cL/t\cL$ is called
{\em quadratic} or simply {\em reflection equation} (RE) algebra.
This algebra can be defined
for any quasitriangular Hopf algebra $\cH$ and its
representation.
When $\cH$ is a quantized
universal enveloping algebra of an algebraic matrix group $G$,
then a certain quotient of the quadratic RE algebra yields a
(one-parameter) quantization of $\kk[G]$.
The mRE algebra $\cL$ as a two-parameter quantization of the
coordinate ring on the matrix space is special for the case $\g=\gl$.
\end{rem}

Let us describe the center $\Inv(\cL)$ of the algebra $\cL$.
First of all, $\Inv(\cL)$  coincides with
the subalgebra of $\Ughr$-invariants and it is isomorphic to
$\Inv_0\ot \kk\dl\hb\dr[t]$, where $\Inv_0\subset\kk[\End(\kk^n)]$ is the
subalgebra of classical invariants.
To describe $\Inv(\cL)$, consider the
matrix $R^*:=\bigl((R^{t_1})^{-1}\bigr)^{t_1}=R^*_1\ot R^*_2$,
where $t_1$ means transposition of the first tensor
component. This matrix is equal to $\cR_1\ot\ga(\cR_2)$ evaluated
in the basic representation.
Define the matrix
$D:=\nu(R_1^*R_2^*)\in \End(\kk^n)\dl\hb\dr$,
where $\nu$ is a scalar.
It is convenient to choose $\nu$ such that
$\dsply\tr(D)=\frac{1-q^{-2n}}{1-q^{-2}}$.
Put $\tau_m=\tr_q(\gen^m):=\tr(D\gen^m)\in \Inv(\cL)$ for $m=1,2,\dots$.
Then the $\Inv(\cL)$ is a polynomial $\K$-algebra
generated by $\tau_1,\dots,\tau_{n-1}$.

\subsection{mRE algebras and twist}

The quantum group $\Ughr$ is a twist of the standard quantization
$\Ugh$.
Let $\cF\in \Ugh^{\ot 2}$ be the corresponding twisting cocycle.
It is an invertible element satisfying the identities
\begin{equation}\label{eq:cocycle}
(\Delta\ot \id)(\cF)\cF_{12}=(\id \ot \Delta)(\cF)\cF_{23},
\end{equation}
\begin{equation}
\label{eq:cocycle_norm}
(\varepsilon\ot \id)(\cF)=1\ot 1=(\id \ot \varepsilon)(\cF),
\end{equation}
where $\varepsilon$ is the counit in $\Ugh$.
As an associative algebra, $\Ughr$ coincides with $\Ugh$
but has a different comultiplication,  $x\mapsto\cF^{-1} \Delta(x)\cF$.
The antipode is transformed accordingly, see \cite{DrinfeldW:AlmostCocomHopfAlgs}.

Recall that
a twist of Hopf algebras induces a transformation of module
algebras, which we also call twist.
Given a $\Ugh$-algebra $\cA$
one gets an algebra over $\Ughr$ with the new multiplication
$a\ot b\mapsto(\cF_{1} a)(\cF_{2} b)$ for $a,b\in \cA$.

Applying $\cF$ to the mRE
algebra, one apparently destroys the form of relations (\ref{eq:mRE}).
Nevertheless, the mRE algebras corresponding to
$\Ughr$ and $\Ugh$
are still related by $\cF$, as we now demonstrate.

Let $\cA$ and  $\wtl \cA$ denote the quadratic RE algebras
corresponding to quantum groups $\Ugh$ and  $\Ughr$ respectively.
We assume that they are extended trivially to $\kk[t]$-algebras.
Let $\wtl \cA'$ be the twist of the algebra $\cA$ by
the cocycle  $\cF$.

Denote by $\{Q_{ij}\}\subset \cA$ and $\{\wtl Q_{ij}\}\subset \wtl\cA$
the generators satisfying the quadratic RE, i.e. the equation similar to
(\ref{eq:mRE}) but with zero in the r.h.s. These algebras
are quantizations of  $\kk[\End(\kk^n)]$ along the STS brackets. As was shown in
\cite{KulishPMudrovA:dRE}, $\wtl\cA'$ is isomorphic to $\wtl\cA$
as $\Ughr$-module algebras,
and the isomorphism $\phi\colon \wtl\cA\to\wtl\cA'$ is given by the
formula
\begin{equation}\label{eq:RETwist}
\widetilde{Q}\mapsto (\rho\circ\ga)(\cF_{1} \ze)Q\rho(\cF_{2}).
\end{equation}
Here $\ze:=\cF^{-1}_{2}\ga^{-1}(\cF^{-1}_{1})\in \Ugh$ is the
element which  participates in definition of the antipode $\tilde{\ga}$ of $\Ughr$,
namely $\tilde \ga(x)=\ga(\ze^{-1} x\ze)$
for all $x\in \Ughr$.

Choose the new generators $\{K_{ij}\}\subset \wtl\cA'$ by setting
$K_{ij}:=Q_{ij}-t\de_{ij}$ and similarly for $\{\wtl K_{ij}\}\subset \wtl\cA$.
Note that $K_{ij}$ are also generators of $\cA$.
\begin{lm}\label{lm:LinGen}
The isomorphism $\phi$ given by the formula (\ref{eq:RETwist}) defines a linear map
$\Span\bigl(\wtl K_{ij}\bigr)\to \Span(K_{ij})$ through the formula
\begin{equation}
\label{Phi}
(\id\ot \phi)(\wtl K)= (\Phi\ot \id)(K),
\end{equation}
where $\Phi$ is an invertible linear operator $\End(V)\to \End(V)$ acting by the rule
$\Phi(X)=(\rho\circ\ga)(\cF_{1}\ze)X\rho(\cF_{2})$
\end{lm}
\begin{proof}
Evaluating $\phi$ on the generators we find
$$
\phi(\wtl K)=\rho\big(\ga(\cF_{1}\ze)\big)Q\rho(\cF_{2})-t=
\rho\big(\ga(\cF_{1}\ze)\big)K\rho\big(\cF_{2}\big)
+t\rho\big(\ga(\cF_{1}\ze)\cF_{2}\big)-t.
$$
The assertion will be proved if we show that $\ga(\ze)\ga(\cF_1)\cF_2=1$.
But this is a well known fact from the twist theory, see \cite{DrinfeldW:AlmostCocomHopfAlgs}.
\end{proof}

Denote by $\cL$ and by $\wtl \cL$ the mRE algebras corresponding to
$\Ugh$ and $\Ughr$, respectively.
\begin{prop}\label{prop:REtwisted}
The algebra $\wtl \cL$ is isomorphic to the twist of 
$\cL$ by the cocycle $\cF$.
\end{prop}
\begin{proof}
Let $\wtl\cA'$ and $\wtl \cL'$ be respectively the twists of the
algebras $\cA$ and $\cL$ by the cocycle $\cF$.
The algebra $\cA$
admits an embedding in $\cL$ through the assignment
\begin{equation}
\label{embed}
 K\mapsto (1-q^{-2}) \gen.
\end{equation}
This embedding induces an embedding $\wtl\cA'\hookrightarrow
\wtl\cL'$ of the  twisted algebras. Let us prove that the isomorphism
(\ref{eq:RETwist}) extends to an isomorphism $\wtl \cL\to\wtl \cL'$.

Denote by $\kk(\!(\hbar)\!)$ the field of Laurent formal series in
$\hbar$. First of all notice that the mapping (\ref{embed}) is
invertible over $\kk(\!(\hbar)\!)$.
Further, the mapping (\ref{eq:RETwist})  induces the isomorphism
$$
\wtl \cL\ot_{\K}\kk(\!(\hbar)\!)\simeq
\wtl\cA\ot_{\K}\kk(\!(\hbar)\!) \longrightarrow
 \wtl\cA'\ot_{\K}\kk(\!(\hbar)\!)\simeq
\wtl\cL'\ot_{\K}\kk(\!(\hbar)\!).$$ of $\kk(\!(\hbar)\!)$-algebras,
which we denote by $\hat \phi$.
Since $\wtl\cL$ and $\wtl\cL'$ are
free over $\kk\dl\hb\dr$, we have the inclusions
$\wtl \cL\subset \wtl \cL\ot_{\K}\kk(\!(\hbar)\!)$
and $ \wtl\cL'\subset \wtl\cL'\ot_{\K}\kk(\!(\hbar)\!)$.
It is therefore sufficient to check that the
image of $\wtl \cL$ under
$\hat \phi$ lies in $ \wtl\cL'$ and similarly for the inverse of $\hat \phi$.

Introduce the linear operator $\Phi^*\colon\Span(K_{ij})\to \Span(K_{ij})$
through the equality $(\Phi\ot \id)(K)=(\id\ot \Phi^*)(K)$
(the dual conjugate of $\Phi$).
Evaluate $\hat \phi$ on a monomial in the generators $\wtl L_{ij}$:
$$
\hat\phi (\wtl L_{i_1 j_1}\ldots \wtl L_{i_k j_k})
=
\hat\phi \left(\frac{1}{\omega}\wtl K_{i_1 j_1}\ldots \frac{1}{\omega}\wtl K_{i_k j_k}\right)
=
\frac{1}{\omega}\Phi^*(K_{i_1 j_1})\ldots \frac{1}{\omega}\Phi^*(K_{i_k j_k}),
$$
where $\omega=1-q^{-2}$.
The last equality is obtained using Lemma~\ref{lm:LinGen}.
But the rightmost expression is
$\Phi^*(L_{i_1 j_1})\ldots \Phi^*(L_{i_k j_k})\in \wtl \cL'$.
In the same fashion, one can check that ${\hat \phi}^{-1}(\wtl \cL') \subset \wtl \cL$.
\end{proof}

\begin{cor}\label{cor:twist}
Let $\ola$ be a semisimple coadjoint orbit. For any quantum group $\Ughr$ there
exists a $\Ughr$-equivariant quantization of $\fun$ which is a
quotient of the mRE algebra associated with $\Ughr$.
\end{cor}
\begin{proof}
Let $\cB$ be the quantization of $\fun$ corresponding to the
standard quantum group $\Ugh$.
It is a quotient of the mRE algebra
$\cL$. The twisted module algebra $\wtl{\cB'}$ is a
$\Ughr$-quantization of $\fun$.
It is a quotient of the algebra
$\wtl\cL'$, which is isomorphic to $\wtl\cL$ by
Proposition~\ref{prop:REtwisted}.
\end{proof}

\subsection{More on RE algebras and twists}
\label{subsect:Twist&REA}

We are going to derive a description of quantum orbits for
an arbitrary quantum group from that corresponding to
the Drinfeld-Jimbo quantum group.
To this end, we need
some facts about Hopf algebras.

As we argued in the previous section (see Proposition~\ref{prop:REtwisted}),
the twist of the (modified) reflection
equation algebra associated with a quantum group
is isomorphic to the (modified) reflection equation algebra
associated with the twisted quantum group.
In this section we obtain a more detailed information about
that isomorphism.
We start with the following auxiliary algebraic assertion.
\begin{lm}\label{4-3}
Let $\cH$ be a Hopf algebra with multiplication $m$, comultiplication $\Delta$,
and invertible antipode $\ga$. Suppose $\cF\in \cH\ot \cH$ is a twisting cocycle.
Then
$$
m_{23}\circ\ga_3\Bigl((\Delta\ot \Delta)(\cF)(\cF\ot\cF)(\ze\ot 1\ot \ze\ot 1)\Bigr)=
\cF_{1}\ze\ot 1\ot \cF_{2},
$$
where the argument in the left-hand-side belongs to $\cH^{\ot 4}$.
\end{lm}
\begin{proof}
Applying the cocycle  equation (\ref{eq:cocycle}) to
$(\Delta\ot \Delta)(\cF)\cF_{34}$, we obtain for the left-hand side
the expression
$$
\cF_1^{(1)}\cF_{1'}^{(1)}\cF_{1''}\ze\ot \cF_1^{(2)}\cF_{1'}^{(2)}\cF_{2''}\ga(\cF_1^{(3)}\cF_{2'}\ze)\ot \cF_2.
$$
In order to distinguish between different copies of $\cF$, the subscripts are marked with dashes.
We apply the cocycle equation to $\cF_{1'}^{(1)}\cF_{1''}\ot \cF_{1'}^{(2)}\cF_{2''}\ot \cF_{2'}$
and obtain
$$
\cF_1^{(1)}\cF_{1'}\ze\ot \cF_1^{(2)}\cF_{2'}^{(1)}\cF_{1''}\ga(\cF_1^{(3)}\cF_{2'}^{(2)}\cF_{2''}\ze)\ot \cF_2.
$$
Now the  statement immediately follows from the equalities
$
\cF_{1}\ga(\ze)\ga(\cF_{2})=1
$
and (\ref{eq:cocycle_norm}).
\end{proof}

Suppose that $\cH$ is a quasitriangular Hopf algebra and
let $(V,\rho)$ be a finite dimensional  representation of $\cH$.
We say that a matrix  $A\in \End(V)\ot \cA$ is {\em invariant}, if
$h\triangleright A=\rho\bigl(\ga(h^{(1)})\bigr)A \rho(h^{(2)})$ for
all $h\in \cH$, where $h\triangleright A$ denotes the action (\ref{eq:action}).
Let $\cA$ and $\wtl \cA$ be the (quadratic) RE algebras
corresponding to
the  Hopf algebras $\cH$ and $\wtl \cH$, where $\wtl \cH$
is the twist of $\cH$ by the cocycle  $\cF$.
The map (\ref{eq:RETwist}) implements an equivariant isomorphism
of $\wtl\cH$-module algebras $\wtl \cA\to \wtl \cA'$ where $\wtl \cA'$
is the twist of $\cA$ by $\cF$.
We can also consider $\phi$ as an isomorphism $\wtl \cA\to \cA$ of $\cH$-modules.

For an invariant matrix $\wtl A\in \End(V)\ot\wtl \cA$
we have
\begin{equation}\label{eq:AtA}
(\id \ot \phi)(\wtl A)=(\Phi\ot \id)(A),
\end{equation}
 where
$A$ is an
invariant matrix in $\End(V)\ot\cA$.
(For the definition of the operator $\Phi$, see  Lemma \ref{lm:LinGen}).
\begin{prop}\label{A*B}
Suppose that  $\wtl A$ and $\wtl B$ are invariant matrices from  $\End(V)\ot\wtl \cA$.
Then
 $(\id\ot \phi)(\wtl A\wtl B)=(\Phi\ot\id) (AB)$,
where $A$ and $B$ are invariant matrices from $\End(V)\ot\cA$ defined by
(\ref{eq:AtA}).
\end {prop}
\begin{proof}
Follows from Lemma~\ref{4-3}.
\end{proof}

For any invariant matrix $A\in \End(V)\ot\cA$
we define an invariant (hence central) element
$\tr_q(A):=\tr_V\bigl(\cR_1\ga(\cR_2)A\bigr)\in \cA$.
Note that  here we suppress
the representation symbol and we do not care
about the normalizing scalar, contrary to  Section~\ref{subsection:mREA}.
\begin{prop}\label{tr-tr0}
Suppose that $\wtl A$ and $A$ are invariant matrices with coefficients
in $\wtl\cA$ and $\cA$, respectively, related by (\ref{eq:AtA}). Then
$$
(\tr_q\ot \phi)(\wtl A)=\tr_q(A).$$
\end{prop}
\begin{proof}
Suppressing the representation symbol $\rho$, we find
$$
(\tr_q\ot \phi)(\wtl A)=\tr\Bigl(\wtl \cR_1\wtl\ga(\wtl \cR_2)\ga(\cF_1\ze)A\cF_2\Bigr),
$$
where $\wtl \ga$ is the antipode in $\wtl \cH$, $\wtl \ga(x)=\ga(\ze^{-1}x\ze)$.
But
$$
\cF_2 \wtl\cR_1\wtl \ga(\wtl \cR_2)\ga(\cF_1\ze)=
\cF_2\wtl\cR_1\ga(\cF_1\wtl \cR_2\ze)
=\cR_1\cF_1\ga(\cR_2\cF_2\ze)
=\cR_1\ga(\cR_2)
$$
because of the equality $\cF_1\ga(\ze)\ga(\cF_2)=1$.
This proves the assertion.
\end{proof}

Denote by $\{Q_{ij}\}\subset \cA$ the RE generators
considered simultaneously as generators for $\wtl \cA'$
(the latter coincides with  $\cA$ as an $\cH$-module
and has the same system of generators as an algebra).
Let $\{\wtl Q_{ij}\}$ denote the RE generators of $ \wtl \cA$.
The matrices $Q$ and $\wtl Q$ are invariant and so are
their powers relative to the multiplications in
$\cA$ and $\wtl \cA$, respectively.
The isomorphism $\phi$ relates $\wtl Q$ and $Q$
by the formula (\ref{eq:AtA}).
The following result is an immediate corollary of Propositions~\ref{A*B} and \ref{tr-tr0}.
\begin{prop}\label{tr-tr}
Regard the algebra isomorphism $\phi\colon \wtl\cA\to \wtl\cA'$ as an isomorphism $\wtl \cA\to \cA$
of vector spaces. Then
$
(\tr_q\ot \phi)(\wtl Q^m)=\tr_q(Q^m).$
\end{prop}

Now let $\wtl \cL$ and $\cL$ be the mRE algebras corresponding
to $\wtl \cH$ and $\cH$.
Let $\{\wtl L_{ij}\}\subset \wtl \cL$ and $\{L_{ij}\}\subset \cL$
be their mRE generators. Put  $\wtl \cL'$ to be  the twist of $\cL$ by $\cF$.
Regard the algebra isomorphism $\phi\colon\wtl \cL\to \wtl \cL'$
extending the isomorphism $\phi\colon\wtl \cA\to \wtl \cA'$
as an isomorphism $\wtl \cL\to \cL$ of vector spaces.
\begin{prop}\label{phi_preserves_submodules}
\begin{enumerate}
\item[{\em (a)}]
The map $\phi$ preserves q-traces:
$
(\tr_q\ot \phi)(\wtl L^m)=\tr_q(L^m).
$
\item[{\em (b)}]
For any polynomial $P$ in one variable,
$(\id \ot \phi)\bigl(P(\wtl L)\bigr)=(\Phi \ot \id)\bigl(P(L)\bigr)$.
\end{enumerate}
\end{prop}
\begin{proof}
The proof readily follows from
Propositions \ref{A*B} and \ref{tr-tr} and
the fact that the twist extends from the quadratic RE algebras to
the modified RE algebras, by Proposition \ref{prop:REtwisted}.
\end{proof}

\section{$\Ughr$-equivariant quantization of orbits}
\label{sect:DMQuantization}

In this section we give a description of a 2-parameter
quantization of the function algebra $\fun$ starting from an
arbitrary (factorizable) classical r-matrix.
This generalizes the construction
given in \cite{DoninJMudrovAI:ExplicitQuantiz}.
The linear term of
this quantization (or, more precisely, the "quasi-Poisson part" of
it (see formula (\ref{eq:UhEquivariantPb})) is
$\hb s^0_\la+t\ka_\la$ where $\hb$ and $t$ are formal parameters.
Reducing this to a one-parameter quantization corresponding to the
curve $t=\la_1\left(e^{-2\hb}-1\right)$ on the plane $(\hb,t)$, we
get a quantization $\funh$ with the linear term
$\hb\left(s^0_\la-2\la_1\ka_\la\right)$.

\subsection{Algebraic description of coadjoint orbits}

Organize the generators of the symmetric algebra $\Sg$ in an
$n\times n$ matrix $\gen=(\gen_{ij})$, then $\Sg=\kk[\gen_{ij}]$.
The algebra $\fun$ of polynomial functions on $\ola$ is a quotient
of $\kk[\gen_{ij}]$ by two sets of relations. The first set of $n^2$
relations can be written in the matrix form as
\begin{equation}\label{eq:CH}
\left(\gen-\la_1\right)\dots\left(\gen-\la_l\right)=0,
\end{equation}
where $(x-\la_1)\dots(x-\la_l)$ is the minimal polynomial for $\la$.
To distinguish the orbits corresponding to the same eigenvalues with
different multiplicities, one should impose the following {\dfn trace
conditions}:
\begin{equation}\label{eq:mults}
\tr\left(\gen^r\right)=\sum_{j=1}^l n_j\la_j^r,\quad r=1,\dots,l-1,
\end{equation}
where $\dsply\sum_{j=1}^l n_j=n$. It is known that the ideal generated
by (\ref{eq:CH}) and (\ref{eq:mults}) is radical, hence it is
precisely the ideal of functions vanishing on $\ola$.

\subsection{On central characters of the mRE algebra}
\label{sec:DMquantization}

To describe quantum orbits explicitly, we need $q$-analogs of the
polynomials in the right-hand side of (\ref{eq:mults}), i.e. quantum
trace functions. For every $m\in \mathbb{N}$ put
$\dsply\hat m:=\frac{1-q^{-2n}}{1-q^{-2}}$.
Fix
$\bla:=(\la_1,\ldots,\la_l)$ and $\bnh:= (\hat n_1,\ldots,\hat n_l)$
assuming $\la_i$ pairwise distinct; put also
$\tl\bla=(\tl\la_1,\ldots,\tl\la_l)$, where
$\dsply\wtl\la_i=\la_i-\frac{t}{\vpi}$.
Consider the family of
functions $\vt_r(\bla,\bnh,q^{-2},t)$, $r=0,\ldots,\infty$, defined
by

\begin{equation}\label{eq:qmults}
\vt_r(\bla, \bnh,q^{-2},t):=\sum_{i=1}^l C_i(\wtl\bla ,
\bnh,\vpi)\la_i^r,
\end{equation}
where
\begin{equation}
C_i(\bla , \bnh,\omega):=\hat n_i\prod_{j\not =i}
\Bigl(1+\vpi\frac{\hat n_j \la_j}{\la_i-\la_j}\Bigr).
\end{equation}
(Recall that we use the notation $\vpi=1-q^{-2}$).
Although manifestly rational, the functions $\vt_r$ are in fact
polynomials in all arguments, see
\cite{DoninJMudrovAI:ExplicitQuantiz}.
In the classical limit $\omega\to 0$, the
function $\vt_r$ turns into the classical trace function
$\dsply\sum_{i=1}^l n_i\la_i^r$.

Fix a polynomial  $\Prj$  in one variable with coefficients in $\kk$.
Consider the quotient of $\cL$ by the
$\Ughr$-invariant ideal of relations $\Prj(\gen)=0$.
Denote by $\Inv_\Prj$ its subalgebra of invariants.
We call  a homomorphism $\Inv_\Prj\to\kk\dl\hb\dr[t]$
{\em a character} of $\Inv_\Prj$.
The meaning of the functions $\vartheta_r(\bla, \bnh,q^{-2},t)$
is explained by the following proposition.
\begin{prop}\label{center_in_M_p}
The algebra $\Inv_P$ is a free module over $\kk\dl\hb\dr[t]$.
The characters of $\Inv_\Prj$ are given by the formulas
\begin{equation}\label{eq:char}
\chi_{\hat \bn}\colon\tau_r\mapsto\vartheta_r(\bla, \bnh,q^{-2},t),
\quad r=1,\ldots, \infty,
\end{equation}
and define an embedding of $\Inv_\Prj$ in the direct sum
$\bigoplus_{\bnh}\kk\dl\hb\dr[t]$.
This embedding becomes an isomorphism
over $\kk\dl \hbar,t\dr$.
\end{prop}

This proposition is proved in \cite{DoninJMudrovAI:ExplicitQuantiz}
for the case of standard quantum group.
One can prove it for
the general quantum group $\Ughr$ using similar arguments.

\subsection{The DM quantization of coadjoint orbits}

Now we are in possession of all ingredients for construction of
quantum orbits. We will work over the ring of scalars being $\kk\dl \hbar,t\dr$.
\begin{thm}\label{th:DMQuantization}
Let $\Ughr$ be {\em any} quasitriangular quantization of $\Ug$ along
a factorizable Lie bialgebra $\g$. Let
let $\kk_{\hbar,t}[\gd]$ be the corresponding mRE algebra generated
by $n^2$ entries of the matrix $\gen$. Then the quotient of
$\kk_{\hbar,t}[\gd]$ by the ideal of relations
\begin{equation}\label{eq:MinPolynomCondition}
(\gen-\la_1)\ldots(\gen-\la_l)=0,
\end{equation}
\begin{equation}\label{eq:qTraceCondition}
\tr_q(\gen^r)=\vartheta_r(\bla, \bnh, q^{-2},t,), \quad
m=1,\ldots,l-1,
\end{equation}
is a $\Ughr$-equivariant quantization of the orbit of matrices with
eigenvalues $\bla$ of multiplicities $\bn$.
\end{thm}
\begin{proof}
The description of the quantized ideal of the orbit can be deduced
from Corollary \ref{cor:twist} and Proposition \ref{center_in_M_p} using
deformation arguments.
We will give an alternative proof based on the results
of Section~\ref{subsect:Twist&REA}, deriving the quantized ideal
of the orbit from the Drinfeld-Jimbo case.

Let $\wtl \cL$ and $\cL$ denote the mRE algebras corresponding
to $\Ughr$ and $\Ugh$, respectively. The quantum group
$\Ughr$ is the twist of $\Ugh$ by a cocycle $\cF$. Denote by $\wtl \cL'$
the corresponding twist of $\cL$; that is a module algebra over $\Ughr$.
By Proposition \ref{prop:REtwisted}, there is an  equivariant
isomorphism of algebras $\phi\colon \wtl \cL\to \wtl \cL'$.
The map $\phi$ is determined by formula (\ref{Phi}), where
the matrices $\wtl K$ and $K$ should be replaced by
$\wtl L$ and $L$, respectively.

Denote by  $\cB$ the quantization of the orbit $O_\la$ which is
equivariant under $\Ugh$.
It is a quotient of $\cL$ by the ideal $\cJ$ of relations
(\ref{eq:MinPolynomCondition})
and (\ref{eq:qTraceCondition}).
The twist $\wtl \cB'$ of the algebra $\cB$ by $\cF$ is a quantization
of $O_\la$ which is equivariant under $\Ughr$.
It is a quotient of $\wtl \cL'$ by the ideal $\wtl \cJ'$ which
coincides with $\cJ$ as a vector space.
Moreover, $\wtl \cJ'$ is generated by the same submodule as $\cJ$ in $\cL$.
In our case that submodule is  spanned by the elements
of the matrix $P(L)$ and the kernel of the central character of $\cL$. Consider the
equivariant isomorphism $\phi^{-1}\colon \wtl \cL'\to \wtl \cL$.
By Proposition \ref{phi_preserves_submodules}, it sends $\Span \bigl(P(L)_{ij}\bigr)$
to $\Span \bigl(P(\wtl L)_{ij}\bigr)$ and preserves the q-traces.
This proves the theorem.
\end{proof}

\begin{rem}
In \cite{MudrovA:QConjClasses}, a description similar to
Theorem \ref{th:DMQuantization} of semisimple quantum conjugacy classes
of the Drinfeld-Jimbo matrix quantum groups is given.
Using the same arguments as in the proof of Theorem \ref{th:DMQuantization}
and the results of  Section \ref{subsect:Twist&REA},
the quantization of \cite{MudrovA:QConjClasses} extends to arbitrary
quantum groups of the classical series.
\end{rem}

\section{Quantization of orbit bundles in $\gld$}
\label{sect:QnOfBundles}

In this section we prove that all orbit bundles admit
$\Ughr$-equivariant quantization and give the explicit construction.
We start with the following
algebraic lemma \cite{DoninJMudrovAI:MethodQuantChar} which we prove
here for the sake of completeness.
\begin{lm}\label{lm:AlgLemma}
Let $Q(x)$ be a polynomial over a field $F$ of zero characteristic,
$\al, \be$ some elements of $F$, and $L,S$ elements of an
associative algebra with unit over $F$ satisfying the following
conditions:
\begin{enumerate}
\item[\em (a)]
$[SLS,L]=0$,
\item[\em (b)]
$S^2=\al S+1$,
\item[\em (c)]
$LQ(L)=\be L$.
\end{enumerate}
Then one has $[SQ(L)S,Q(L)]=0$.
\end{lm}
\begin{rem}
The algebra generated by $S$ and $L$ subject to conditions (a)--(c) is a special case of  cyclotomic affine Hecke algebra of rank 1.
\end{rem}
\begin{proof}
Prove, using the induction on $m\geq 1$, that $[SL^mS,Q(L)]=0$. The
induction base, $m=1$, holds true for one checks readily that (a)
implies $[SLS,L^k]=0$ for any $k$. Now, suppose $[SL^mS,Q(L)]=0$,
then one has using (b):
\begin{gather}
[SL^{m+1}S,Q(L)]=[SL1L^mS,Q(L)]=[SL(S^2-\al S)L^mS,Q(L)]=\nonumber \\
=[SLS^2L^m S,Q(L)]-\al[SLSL^mS,Q(L)].\label{eq:One}
\end{gather}
According to the induction assumption, both $SLS$ and $SL^mS$ commute
with $Q(L)$, thus
$$
[SLS^2L^mS,Q(L)]=[(SLS)(SL^mS),Q(L)]=0.
$$
The last term in (\ref{eq:One}) is treated as follows:
\begin{eqnarray*}
[SLSL^mS,Q(L)]&=&
SL\bigl(SL^mSQ(L)\bigr)-\bigl(Q(L)SLS\bigr)L^mS\\
&=&SLQ(L)SL^mS-SLSL^mQ(L)S \\
&=&\be SLSL^mS-\be SLSL^mS=0,
\end{eqnarray*}
where the induction assumption and (c) were used.
\end{proof}

Recall from see Proposition~\ref{prop:P} that any orbit map is determined
by a polynomial $\Prj$ in one variable.
\begin{thm}\label{th:Main}
Fix a factorizable quantum group $\Ughr$, where $\g=\gl$.
Let $\ola$
and $\omu$ be two orbits in $\g$ satisfying the condintions of
Theorem~\ref{th:MainPoiss}, and denote by
$\kk_\hb[\ola]=\kk_\hb[\gen_{\la,ij}]$ and
$\funh[\omu]=\kk_\hb[\gen_{\mu,ij}]$ their quantizations from
Theorem~\ref{th:DMQuantization} with
$t=\la_1\left(e^{-2\hb}-1\right)$.
Then the assignment
$\gen_\mu\mapsto \Prj(\gen_\la)$, where the polynomial $\Prj$ is
given by (\ref{eq:DefOfP}), is a $\Ughr$-equivariant quantization of
the orbit bundle $\ola\to\omu$ determined by $P$.
\end{thm}
\begin{proof}
Denote by $\prj^*$ the algebra monomorphism $\fun[\omu]\to\fun$
corresponding to the map $\prj$.
Both $\gen_\la$ and $\gen_\mu$ are
subject to the relations (\ref{eq:CH}) and (\ref{eq:mults}).
The algebra homomorphism $\prj^*\colon \fun[\omu]\to\fun$ is determined
by the correspondence $\gen_\mu\mapsto P(\gen_\la)$.
We need to prove that
the same correspondence defines a $\kk\dl\hb\dr$-algebra monomorphism
$\funh[\omu]\to\funh$, i.e. that the matrix $\prj(\gen_\la)$
satisfies the same relations as the matrix $\gen_\mu$.

1. Check the relation: $[SP(\gen_\la)
S,P(\gen_\la)]=\mu_1\hcke[S,P(\gen_\la)]$.

It can be written in the form:
\begin{equation}\label{eq:REforP}
[S\left(P(\gen_\la)-\mu_1\right) S,P(\gen_\la)-\mu_1]=0,
\end{equation}
as $S$ is a Hecke matrix.
It is easy to check that $(\gen_\la-\la_1)\left( P(\gen_\la)-\mu_1
\right)=(\mu_2-\mu_1)(\gen_\la-\la_1)$. Now set $\be:=\mu_2-\mu_1$,
$L:=\gen_\la-\la_1$, $Q(x):=P(x+\la_1)-\mu_1$, then
(\ref{eq:REforP}) follows from Lemma~\ref{lm:AlgLemma}.

2. Check the relation:
\begin{equation}\label{eq:CHforP}
(P(\gen_\la)-\mu_1)(P(\gen_\la)-\mu_2)=0.
\end{equation}
Substituting (\ref{eq:DefOfP}) into the l.h.s. of (\ref{eq:CHforP}),
one gets
\begin{equation}\label{eq:CheckingHC}
\prod_{i=2}^l(\gen_\la-\la_i) \left(
\prod_{i=2}^l(\gen_\la-\la_i)-\prod_{i=2}^l(\la_1-\la_i) \right)
\end{equation}
up to a constant multiple. The expression in the big brackets is
divisible by $\gen_\la-\la_1$.
Indeed, for any polynomial $f(x)$,
the polynomial in two variables  $F(x,y):=f(x)-f(y)$ is divisible by
$x-y$.
This implies that (\ref{eq:CheckingHC}) is divisible by the
minimal polynomial of $\la$, so it is equal to zero.

3. In order to check the q-Trace Condition,
\begin{equation}\label{eq:TraceCondForP}
\tr_q P(\gen_\la)=\tr_q(\gen_\mu),
\end{equation}
we put $\bnu=\bnh=(\hat n_1,\ldots,\hat n_l)$ and $\vpi:=1-q^{-2}$
in the functions $C_i(\bla,\bnu,\vpi)$, see Appendix, formula
(\ref{eq:Cs}). (As above, $\hat n_i =
\frac{1-q^{-2n_i}}{1-q^{-2}}$).

Replacing  $\bla$, $\bmu$, $\gen_\la$, $\gen_\mu$ and $\Prj(x)$
by $(0, \la_2-\la_1,\ldots,\la_l-\la_1)$, $(0,\mu_2-\mu_1)$,
$\gen_\la-\la_1$, $\gen_\mu-\mu_1$ and $\Prj(x+\la_1)-\mu_1$
respectively, one reduces the problem to the case $\la_1=\mu_1=0$.
So, it suffices to prove that the condition (\ref{eq:TraceCondForP})
is satisfied when $\la_1=0$ and $\Prj(0)=0$.

By assumption, $\mu_1=\Prj(0)=0$, therefore one has:
\begin{equation}\label{eq:TrMu}
\qtr \gen_\mu=C_2(\bmu,\bmh,q)\mu_2=\hat n'\mu_2,
\end{equation}
where $\bm=(n_1,n')$ and $\dsply n':=\sum_{i=2}^l n_i$. On the other
hand,
\begin{equation}\label{eq:TrLa}
\qtr\bigl(P(\gen_\la)\bigr)= \sum_{i=1}^l P(\la_i)C_i(\bla,\bnh,q)=
\mu_2\sum_{i=2}^l C_j(\bla,\bnh,q),
\end{equation}
because $\Prj(\la_1)=\Prj(0)=0$. By Corollary~\ref{prop:Reductio}
(see Appendix),
\begin{equation}\label{eq:SumOfCs}
\dsply\sum_{i=2}^l C_i(\bla, \bnh, \vpi)=\hat n',
\end{equation}
since $1-\vpi\hat n_i=q^{-2n_i}$. Substituting (\ref{eq:SumOfCs})
into (\ref{eq:TrLa}), one concludes that the latter is equal to
(\ref{eq:TrMu}).
\end{proof}

\section*{Appendix}

In this section, we study some properties of the coefficients $C_i$
in (\ref{eq:qmults}), which were announced without proof in
\cite{DoninJMudrovAI:GeneralizedVermaMod}.
For $1\leqs i\leqs l$,
define a function of $2l+1$ variables $\bla=(\la_1,\ldots,\la_l)$,
$\bnu=(\nu_1,\ldots,\nu_l)$ and $\vpi$:
\begin{equation}\label{eq:Cs}
C_i(\bla,\bnu,\vpi):= \nu_i\prod_{\substack{1\leqs j\leqs l\\j\not
=i}} \left(1+\vpi\frac{\nu_{j}\la_{j}}{\la_i-\la_j}\right),
\end{equation}
and also another function of the same variables:
\begin{equation*}
S(\bla,\bnu,\vpi)=\sum_{i=1}^l C_i(\bla,\bnu,\vpi).
\end{equation*}
These functions were introduced in
\cite{DoninJMudrovAI:ExplicitQuantiz} and
\cite{DoninJMudrovAI:GeneralizedVermaMod}. Our goal is to prove
Proposition~\ref{prop:MainLemmaOnC} below.

Obviously, $S(\bla,\bnu,\vpi)$ is stable under simultaneous
permutations of the entries of $\bla$ and the  entries of $\bnu$.
In fact, a stronger statement is true:
\begin{lm}
$S(\bla,\bnu,\vpi)$ is a symmetric function of $\bla$.
\end{lm}
\begin{proof}
It suffices to show that $S(\bla,\bnu,\vpi)$ is stable under the
transposition $\la_1\leftrightarrow \la_2$.
First, opening the
brackets in (\ref{eq:Cs}) one gets
\begin{equation}\label{eq:CsExpanded}
C_i(\bla,\bnu,\vpi)= \nu_i+\nu_i\sum_{k=1}^{l-1}\vpi^ k
\sum_{j_1<\ldots<j_ k} \frac{\nu_{j_1}\la_{j_1}}{\la_i-\la_{j_1}}
\ldots \frac{\nu_{j_ k}\la_{j_ k}}{\la_i-\la_{j_ k}}.
\end{equation}
In this form, the functions $C_i$ were introduced in \cite{DoninJMudrovAI:ExplicitQuantiz}.
The multiplicative form (\ref{eq:Cs}) appeared in
\cite{GurewitzSaponov:GeomNCOrbits}.
All the terms in (\ref{eq:CsExpanded}) containing $\la_1$ and
$\la_2$ can be arranged into sums of the following three forms:
\begin{gather*}
\nu_j\frac{\nu_1\la_1}{\la_j-\la_1} \frac{\nu_2\la_2}{\la_j-\la_2}f
=\nu_j\nu_1\nu_2\frac{\la_1}{\la_j-\la_1}
\frac{\la_2}{\la_j-\la_2}f, \\
\ \\
\nu_1\frac{\nu_2\la_2}{\la_2-\la_1}f+
\nu_2\frac{\nu_1\la_1}{\la_1-\la_2}f=
\nu_1\nu_2f, \\
\ \\
\nu_j\frac{\nu_1\la_1}{\la_j-\la_1}f+
 \nu_j\frac{\nu_2\la_2}{\la_j-\la_2}f+
 \nu_1\frac{\nu_j\la_j}{\la_1-\la_j}f+
 \nu_2\frac{\nu_j\la_j}{\la_2-\la_j}f
=-(\nu_j\nu_1+\nu_j\nu_2)f,
\end{gather*}
with $j\not=1,2$, and $f$ being independent on $\la_1$ and $\la_2$.
It is seen that the expressions in the right hand sides are stable
under the transposition $\la_1\leftrightarrow \la_2$.
\end{proof}
\begin{prop}\label{prop:MainLemmaOnC}
$\dsply\vpi S(\bla,\bnu,\vpi)=1-\prod_{i=1}^l(1-\vpi\nu_i)$.
\end{prop}
\begin{proof}
Prove first that $S(\bla,\bnu,\vpi)$ does not actually depend on
$\bla$. Fix $\bnu$ and $\vpi$, and consider $S(\bla,\bnu,\vpi)$ as a
rational function of $\bla$ only. This function is homogeneous of
degree zero.
Reducing $S(\bla,\bnu,\vpi)$ to the common denominator
$\prod_{i<j}(\la_i-\la_j)$ we obtain a ratio of two homogeneous
polynomials of the same degree.
Since $S(\bla,\bnu,\vpi)$ is a
symmetric function of $\bla$, the numerator is divisible by
$\prod_{i<j}(\la_i-\la_j)$ because the ring of polynomials is a
unique factorization domain.
Since the numerator of
$S(\bla,\bnu,\vpi)$  has the same degree as the denominator,
$S(\bla,\bnu,\vpi)$ is independent on $\bla$.

Now put $\la_l=0$, then it follows from (\ref{eq:Cs}) that
$S(\bla,\bnu,\vpi)=S(\bla',\bnu',\vpi)+
\nu_l\prod_{i=1}^{l-1}(1-\vpi\nu_i)$, where
$\bla'=(\la_1,\ldots,\la_{l-1})$ and
$\bnu'=(\nu_1,\ldots,\nu_{l-1})$. Finally, one applies the induction
on $l$.
\end{proof}

\begin{cor}\label{prop:Reductio}
\begin{enumerate}
\item[{\em (a)}]
If $\la_1=0$ then $\dsply\vpi\sum_{i=2}^l C_i(\bla,\bnu,\vpi)
=1-\prod_{i=2}^l (1-\vpi\nu_i)$.

\item[{\em (b)}]
Denote $\bla'=(\la_2,\ldots,\la_l)$,
$\bnu'=(\nu_2,\ldots,\nu_l)\in\kk^{l-1}$, and suppose that
$\la_1=0$.
Then $\dsply\sum_{i=2}^l C_i(\bla,\bnu,\vpi)
=\sum_{i=2}^l C_i(\bla',\bnu',\vpi)$.
\item[{\em (c)}]
One has $\dsply\sum_{i=1}^l C_i(\bla, \bnh, \vpi)=\hat n$.
\end{enumerate}
\end{cor}
\begin{proof}
(a) Denote $S'(\bla,\bnu,\vpi)=\dsply\vpi\sum_{i=2}^l
C_i(\bla,\bnu,\vpi)$.
Then
\begin{gather*}
\vpi S'(\bla,\bnu,\vpi)=
\vpi\bigl(S(\bla,\bnu,\vpi)-C_1(\bla,\bnu,\vpi)\bigr)=\\
=1-\prod_{i=1}^l(1-\vpi\nu_i)-\vpi\nu_1\prod_{i=2}^l
(1-\vpi\nu_i)=1-\prod_{i=2}^l (1-\vpi\nu_i).
\end{gather*}

(b) Obvious.

(c) Note that $1-\vpi\hat n_i=q^{-2n_i}=e^{-2n_i\hb}$, recall that
$\displaystyle n=\sum_{i=1}^l n_i$ and use
Proposition~\ref{prop:MainLemmaOnC}.
\end{proof}

\bibliography{ormath}

\begin{thebibliography}{10}

\bibitem{BelawinDrinfeld:CYBE}
A.~A. Belavin and V.~G. Drinfeld.
\newblock On solutions of {Y}ang---{B}axter equation.
\newblock {\em Functional Analysis and Applications}, 16(3):1--29, 1982.
\newblock {\sl English translation:} {\em Functional Analysis and
  Applications}, vol.32 (1985), p.p. 254--255.

\bibitem{DoninJ:DoubleQuantizationOfLieGroups}
J.~Donin.
\newblock Double quantization on coadjoint representations of simple {L}ie
  groups and its orbits.
\newblock {\em MPIM}, September 1999.

\bibitem{DoninJ:QuanGMans}
J.~Donin.
\newblock Quantum $g$-manifolds.
\newblock {\em Contemp. Math.}, 315:47--60, 2002.

\bibitem{DoninJMudrovAI:ExplicitQuantiz}
J.~Donin and A.~Mudrov.
\newblock Explicit quantization on coadjoint orbits of $gl(n,\mathbb{C})$.
\newblock {\em Lett. Math. Phys}, 62(1):17--32, 2002.

\bibitem{DoninJMudrovAI:MethodQuantChar}
J.~Donin and A.~Mudrov.
\newblock Method of quantum characters in equivariant quantization.
\newblock {\em Comm. Math. Phys.}, 234:533--555, 2003.

\bibitem{DoninJMudrovAI:GeneralizedVermaMod}
J.~Donin and A.~Mudrov.
\newblock Quantum coadjoint orbits of {$GL(n)$} and generalized {V}erma
  modules.
\newblock {\em Lett. Math. Phys.}, 67:167--184, 2004.

\bibitem{DoninJOstapenkoV:EquivQn}
J.~Donin and V.~Ostapenko.
\newblock Equivariant quantization on quotients of simple {L}ie groups by
  reductive subgroups.
\newblock {\em Czech. Journ. of Phys.}, 52(11):1213--1218, 2002.

\bibitem{DoninGurevicShnider:DoubleQuantiz}
J.~Donin{, D. Gurevich} and S.~Shnider.
\newblock Double quantization on some orbits in the coadjoint representations
  of simple {L}ie groups.
\newblock {\em Comm. Math. Phys.}, 204(1):39--60, 1999.
\newblock math.QA/9807159.

\bibitem{DrinfeldW:Doklad}
V.~G. Drinfeld.
\newblock Quantum groups.
\newblock In {\em Proceedings of {I}nternational {C}ongress of
  {M}athematicians, Berkley 1986}, volume~1, pages 798--820. AMS, Providence,
  1986.

\bibitem{DrinfeldW:AlmostCocomHopfAlgs}
V.~G. Drinfeld.
\newblock Almost cocommutative {H}opf algebras.
\newblock {\em Algebra i Analys}, 1(2):321--342, 1989.
\newblock {\sl English translation:} {\em Leningrad Journal of Mathematics},
  vol. 1 No. 2 (1990), p.p. 321--342.

\bibitem{EtinhoffKazdan:QuantBialg}
P.~Etingoff and D.~Kazhdan.
\newblock Quantization of {L}ie bialgebras.
\newblock {\em Selecta Math.}, 2(1):1--41, 1996.

\bibitem{GurewitzSaponov:GeomNCOrbits}
D.~Gurevich and P.~Saponov.
\newblock Geometry of non-commutative orbits related to {H}ecke symmetries.
\newblock math.QA/0411579.

\bibitem{KarolinskijE:ClassPoissHomSpace}
E.~Karolinskii.
\newblock A classification of {P}oisson homogeneous spaces of complex reductive
  {P}oisson--{L}ie groups.
\newblock In P.~Urbanski J.~Grabowski, editor, {\em Poisson geometry},
  volume~51, pages 103--108. Banach Center, Warsaw, 2000.

\bibitem{KulishPMudrovA:dRE}
P.~Kulish and A.~Mudrov.
\newblock Dynamical reflection equation.
\newblock math.QA/0405556.

\bibitem{MudrovA:QConjClasses}
A.~Mudrov.
\newblock Quantum conjugacy classes of simple matrix groups.
\newblock math.QA/0412538.

\bibitem{STS:QDualityQDouble}
M.~Semenov-Tian-Shansky.
\newblock {P}oisson-{L}ie groups, quantum duality principle, and the quantum
  double.
\newblock {\em Contemp. Math.}, 175:219--248, 1994.

\end{thebibliography}

\end{document}